\documentclass[a4paper, 12 pt]{article}

\usepackage[T2A]{fontenc}
\usepackage[cp1251]{inputenc}
\usepackage[english]{babel}
\usepackage[tbtags]{amsmath}
\usepackage{amsfonts,amssymb}

\usepackage{mathrsfs}

\usepackage{geometry} 
\begin{document}
\title{Algebra of formal power series, isomorphic to the algebra of formal Dirichlet series}
 \author{E. Burlachenko}
 \date{}

 \maketitle
\begin{abstract}
Ordinary algebra of formal power series in one variable is convenient to study by means of the  algebra of  Riordan matrices and the Riordan group. In this paper we consider algebra of formal power series without constant term, isomorphic to the algebra of formal Dirichlet series. To study it, we introduce matrices, similar to the Riordan matrices. As a result, some analogies between two algebras becomes visible. For example, the Bell polynomials (polynomials of partitions of number $n$ into $m$ parts) play a certain role in the ordinary algebra. Similar polynomials (polynomials of decompositions of number $n$ into $m$ factors)  play a similar role in the considered algebra. Analog of the Lagrange series for the considered algebra is also exists. In connection with this analogy, we introduce matrix group, similar to the Riordan group and called the Riordan-Dirichlet group. As an example, we consider analog of the Abel's identities for this group.
\end{abstract}
\section{Introduction}
Transformations, corresponding to multiplication and composition of  series,
play the main role in the space of formal power series over the field of real or complex  numbers. Multiplication is geven by the matrix $\left( a\left( x \right),x \right)$ $n$th column of which, $n=0,\text{ }1,\text{ }2,\text{ }...$ ,  has the generating function $b\left( x \right){{x}^{n}}$; composition is given by the matrix $\left( 1,a\left( x \right) \right)$ $n$th column of which  has the generating function ${{a}^{n}}\left( x \right)$, ${{a}_{0}}=0$:
$$\left( b\left( x \right),x \right)g\left( x \right)=b\left( x \right)g\left( x \right), \qquad\left( 1,a\left( x \right) \right)g\left( x \right)=g\left( a\left( x \right) \right).$$
Matrix
$$\left( b\left( x \right),x \right)\left( 1,a\left( x \right) \right)=\left( b\left( x \right),a\left( x \right) \right)$$
is called Riordan array [1] – [4]; $n$th column of Riordan array has the generating function $b\left( x \right){{a}^{n}}\left( x \right)$. Thus,
$$\left( b\left( x \right),a\left( x \right) \right)f\left( x \right){{g}^{n}}\left( x \right)=b\left( x \right)f\left( a\left( x \right) \right){{\left( g\left( a\left( x \right) \right) \right)}^{n}},$$
$$\left( b\left( x \right),a\left( x \right) \right)\left( f\left( x \right),g\left( x \right) \right)=\left( b\left( x \right)f\left( a\left( x \right) \right),g\left( a\left( x \right) \right) \right).$$
Matrices $\left( b\left( x \right),a\left( x \right) \right)$, ${{b}_{0}}\ne 0$, ${{a}_{1}}\ne 0$, form a group called the Riordan group.

$n$th coefficient of the series $a\left( x \right)$, $n$th row and $n$th column of the matrix $A$ will be denoted  respectively by 
$$\left[ {{x}^{n}} \right]a\left( x \right),\qquad   \left[ n,\to  \right]A,   \qquad[\uparrow ,n]A,$$
at that $\left[ {{x}^{n}} \right]a\left( x \right)b\left( x \right)=\left[ {{x}^{n}} \right]\left( a\left( x \right)b\left( x \right) \right)$. We associate rows and columns of matrices with the generating functions of their elements. 

Matrices 
$${{\left| {{e}^{x}} \right|}^{-1}}\left( b\left( x \right),a\left( x \right) \right)\left| {{e}^{x}} \right|={{\left( b\left( x \right),a\left( x \right) \right)}_{{{e}^{x}}}},$$
where $\left| {{e}^{x}} \right|$ is the diagonal matrix whose diagonal elements are equal to the coefficients of  the series ${{e}^{x}}$: $\left| {{e}^{x}} \right|a\left( x \right)=\sum\nolimits_{n=0}^{\infty }{{{{a}_{n}}{{x}^{n}}}/{n!}\;}$, are called exponential Riordan arrays. Denote
$$\left[ n,\to  \right]{{\left( b\left( x \right),a\left( x \right) \right)}_{{{e}^{x}}}}={{s}_{n}}\left( x \right), \qquad{{b}_{0}}\ne 0, \qquad{{a}_{1}}\ne 0.$$
Then
$${{\left( b\left( x \right),a\left( x \right) \right)}_{{{e}^{x}}}}{{\left( 1-\varphi x \right)}^{-1}}={{\left| {{e}^{x}} \right|}^{-1}}\left( b\left( x \right),a\left( x \right) \right){{e}^{\varphi x}}={{\left| {{e}^{x}} \right|}^{-1}}b\left( x \right)\exp \left( \varphi a\left( x \right) \right),$$
or
$$\sum\limits_{n=0}^{\infty }{\frac{{{s}_{n}}\left( \varphi  \right)}{n!}{{x}^{n}}}=b\left( x \right)\exp \left( \varphi a\left( x \right) \right).$$
Sequence of polynomials ${{s}_{n}}\left( x \right)$ is called Sheffer sequence, and in the case $b\left( x \right)=1$, binomial sequence [4]. If ${{s}_{n}}\left( x \right)$ corresponds to the $n$th row of the matrix ${{\left( 1,\log a\left( x \right) \right)}_{{{e}^{x}}}}$, then
$${{a}^{\varphi }}\left( x \right)=\sum\limits_{n=0}^{\infty }{\frac{{{s}_{n}}\left( \varphi  \right)}{n!}{{x}^{n}}}.$$

Matrices 
$$\left( a\left( x \right),x \right)=\sum\limits_{n=0}^{\infty }{{{a}_{n}}}\left( {{x}^{n}},x \right),$$
$$\left( a\left( x \right),x \right)=\left( \begin{matrix}
   {{a}_{0}} & 0 & 0 & 0 & \ldots   \\
   {{a}_{1}} & {{a}_{0}} & 0 & 0 & \ldots   \\
   {{a}_{2}} & {{a}_{1}} & {{a}_{0}} & 0 & \ldots   \\
   {{a}_{3}} & {{a}_{2}} & {{a}_{1}} & {{a}_{0}} & \ldots   \\
   \vdots  & \vdots  & \vdots  & \vdots  & \ddots   \\
\end{matrix} \right)$$
form the algebra, isomorphic to the algebra of formal power series. Theme of this paper appeared as a result of the following observation. If in the algebra of  matrices $\left( a\left( x \right),x \right)=\sum\nolimits_{n=1}^{\infty }{{{a}_{n}}}\left( {{x}^{n}},x \right)$, isomorphic to the algebra of formal power series without constant term, the matrix of multiplication $\left( {{x}^{n}},x \right)$ is replaced by the matrix of composition $\left( 1,{{x}^{n}} \right)$, the result will be algebra, isomorphic to the algebra of formal Dirichlet series.

In the following sections of this paper we will consider the basic elementary aspects of the algebra thus obtained. Emphasis is on its relationship with ordinary algebra of formal power series. Research tools are the matrices, similar to the Riordan matrices. Group, similar to the Riordan group, is introduced in the last section, in which we consider series, similar to the Lagrange series
\section{}
Denote
$\left( 1,{{x}^{n}} \right)=\left\langle {{x}^{n}},x \right\rangle $,
where we  take into consideration the analogy with Riordan matrices, which will be developed in the future. Then
$$\left\langle a\left( x \right),x \right\rangle =\sum\limits_{n=1}^{\infty }{{{a}_{n}}}\left\langle {{x}^{n}},x \right\rangle ,$$
$$\left\langle a\left( x \right),x \right\rangle =\left(\setcounter{MaxMatrixCols}{20} \begin{matrix}
   {{a}_{\Sigma }} & 0 & 0 & 0 & 0 & 0 & 0 & 0 & 0 & 0 & \ldots   \\
   0 & {{a}_{1}} & 0 & 0 & 0 & 0 & 0 & 0 & 0 & 0 & \dots   \\
   0 & {{a}_{2}} & {{a}_{1}} & 0 & 0 & 0 & 0 & 0 & 0 & 0 & \ldots   \\
   0 & {{a}_{3}} & 0 & {{a}_{1}} & 0 & 0 & 0 & 0 & 0 & 0 & \ldots   \\
   0 & {{a}_{4}} & {{a}_{2}} & 0 & {{a}_{1}} & 0 & 0 & 0 & 0 & 0 & \ldots   \\
   0 & {{a}_{5}} & 0 & 0 & 0 & {{a}_{1}} & 0 & 0 & 0 & 0 & \ldots   \\
   0 & {{a}_{6}} & {{a}_{3}} & {{a}_{2}} & 0 & 0 & {{a}_{1}} & 0 & 0 & 0 & \ldots   \\
   0 & {{a}_{7}} & 0 & 0 & 0 & 0 & 0 & {{a}_{1}} & 0 & 0 & \ldots   \\
   0 & {{a}_{8}} & {{a}_{4}} & 0 & {{a}_{2}} & 0 & 0 & 0 & {{a}_{1}} & 0 & \ldots   \\
   0 & {{a}_{9}} & 0 & {{a}_{3}} & 0 & 0 & 0 & 0 & 0 & {{a}_{1}} & \ldots   \\
   \vdots  & \vdots  & \vdots  & \vdots  & \vdots  & \vdots  & \vdots  & \vdots  & \vdots  & \vdots  & \ddots   \\
\end{matrix} \right),$$
where 
$$[\uparrow ,n]\left\langle a\left( x \right),x \right\rangle =a\left( {{x}^{n}} \right),     \qquad a\left( 1 \right)=\sum\limits_{n=1}^{\infty }{{{a}_{n}}}={{a}_{\Sigma }}.$$

Sum of the coefficients of formal power series, obviously, in need of definition. Perhaps, each numerical series $\sum\nolimits_{n=1}^{\infty }{{{a}_{n}}}$ has the value equal to a certain number or $\pm \infty $, which in the case of convergent series coincides with the its  sum. For divergent numerical series, selection of the value corresponding to the sum of series, is ambiguous and depends on the accepted conditions [5]. We go around this problem and will consider the expression ${{a}_{\Sigma }}$ as “formal numerical series”. Actions with the formal numerical series we define by the action with the corresponding power series: if $a\left( x \right)+b\left( x \right)=c\left( x \right)$, then ${{a}_{\Sigma }}+{{b}_{\Sigma }}={{c}_{\Sigma }}$; if  $\left\langle a\left( x \right),x \right\rangle b\left( x \right)=c\left( x \right)$, then ${{a}_{\Sigma }}{{b}_{\Sigma }}={{c}_{\Sigma }}$. Numerical series, corresponding to the series $a\left( {{x}^{n}} \right)$, are considered to be identical; if $a\left( x \right)=\varphi {{x}^{n}}$, then ${{a}_{\Sigma }}=\varphi $.

Algebras of the matrices $\left( a\left( x \right),x \right)$, $\left\langle a\left( x \right),x \right\rangle $ and the corresponding algebras of formal power series we will be called the $\left( a\left( x \right),x \right)$- algebra and the $\left\langle a\left( x \right),x \right\rangle $- algebra. Denote
$$\left\langle a\left( x \right),x \right\rangle b\left( x \right)=a\left( x \right)\circ b\left( x \right), \qquad{{b}_{0}}=0.$$
If $a\left( x \right)\circ b\left( x \right)=c\left( x \right)$, then ${{c}_{n}}=\sum\nolimits_{d|n}{{{a}_{d}}{{b}_{{n}/{d}\;}}}$, where summation is over all divisors $d$ of number $n$. Inverse to the series $a\left( x \right)$ we call  the series ${{a}^{\left( -1 \right)}}\left( x \right)$, which is defined by the identity
$$a\left( x \right)\circ {{a}^{\left( -1 \right)}}\left( x \right)={{a}^{\left( 0 \right)}}\left( x \right)=x.$$
This is consistent with the fact that $\left\langle x,x \right\rangle $ is the identity matrix: ${{x}^{n}}\circ a\left( x \right)=a\left( {{x}^{n}} \right)$. Denote also
$${{a}^{\left( n-1 \right)}}\left( x \right)\circ a\left( x \right)={{a}^{\left( n \right)}}\left( x \right).$$

Note parallels between two algebras. Since
$$\left( {{x}^{n}},x \right)\left( {{x}^{m}},x \right)=\left( {{x}^{m}},x \right)\left( {{x}^{n}},x \right)=\left( {{x}^{n+m}},x \right),$$
$$\left\langle {{x}^{n}},x \right\rangle \left\langle {{x}^{m}},x \right\rangle =\left\langle {{x}^{m}},x \right\rangle \left\langle {{x}^{n}},x \right\rangle =\left\langle {{x}^{nm}},x \right\rangle ,$$
then
$$\left( a\left( x \right),x \right)\left( b\left( x \right),x \right)=\left( b\left( x \right),x \right)\left( a\left( x \right),x \right);$$
$$\left\langle a\left( x \right),x \right\rangle \left\langle b\left( x \right),x \right\rangle =\left\langle b\left( x \right),x \right\rangle \left\langle a\left( x \right),x \right\rangle .$$
Since
$${{\left( {{x}^{m}},x \right)}^{n}}=\left( {{x}^{mn}},x \right),$$
then the matrix $\left( a\left( x \right),x \right)$ is the power series:
$$\left( a\left( x \right),x \right)=\sum\limits_{n=0}^{\infty }{{{a}_{n}}{{\left( x,x \right)}^{n}}};$$
since
$${{\left\langle {{x}^{m}},x \right\rangle }^{n}}=\left\langle {{x}^{{{m}^{n}}}},x \right\rangle ,$$
then
$$\left\langle a\left( x \right),x \right\rangle ={{a}_{1}}\left\langle x,x \right\rangle +\sum\limits_{m=2}^{\infty }{\sum\limits_{n=1}^{\infty }{{{a}_{{{m}^{n}}}}}}{{\left\langle {{x}^{m}},x \right\rangle }^{n}},  \qquad m\ne {{k}^{s}}, \qquad s>1.$$
Thus, matrices of the form
$$\left\langle a\left( x \right),x \right\rangle =\sum\limits_{n=0}^{\infty }{{{a}_{n}}{{\left\langle {{x}^{m}},x \right\rangle }^{n}}}, \qquad m>1,$$
being power series, form the algebra, isomorphic to the $\left( a\left( x \right),x \right)$- algebra: if
$$\left( \sum\limits_{n=0}^{\infty }{{{a}_{n}}{{x}^{n}}} \right)\left( \sum\limits_{n=0}^{\infty }{{{b}_{n}}{{x}^{n}}} \right)=\sum\limits_{n=0}^{\infty }{{{c}_{n}}{{x}^{n}}},$$
then
$$\left( \sum\limits_{n=0}^{\infty }{{{a}_{n}}{{x}^{{{m}^{n}}}}} \right)\circ \left( \sum\limits_{n=0}^{\infty }{{{b}_{n}}{{x}^{{{m}^{n}}}}} \right)=\sum\limits_{n=0}^{\infty }{{{c}_{n}}{{x}^{{{m}^{n}}}}}.$$
For example, for integers $k$,
$${{\left( x+{{x}^{m}} \right)}^{\begin{matrix}
   \left( k \right)  \\
\end{matrix}}}=\sum\limits_{n=0}^{\infty }{\left( \begin{matrix}
   k  \\
   n  \\
\end{matrix} \right){{x}^{{{m}^{n}}}}}.$$

In the $\left( a\left( x \right),x \right)$-algebra the identity
$$\left( \sum\limits_{n=0}^{\infty }{{{a}_{n}}{{\beta }^{n}}{{x}^{n}}} \right)\left( \sum\limits_{n=0}^{\infty }{{{b}_{n}}{{\beta }^{n}}{{x}^{n}}} \right)=\sum\limits_{n=0}^{\infty }{{{c}_{n}}{{\beta }^{n}}{{x}^{n}}}$$
holds for any values $\beta $. In the $\left\langle a\left( x \right),x \right\rangle $-algebra the similar identity holds:
$$\left( \sum\limits_{n=1}^{\infty }{{{a}_{n}}{{n}^{\beta }}{{x}^{n}}} \right)\circ \left( \sum\limits_{n=1}^{\infty }{{{b}_{n}}{{n}^{\beta }}{{x}^{n}}} \right)=\sum\limits_{n=1}^{\infty }{{{c}_{n}}{{n}^{\beta }}{{x}^{n}}}.$$

\section{}

We introduce matrices $\left\langle x|a\left( x \right) \right\rangle $, which will play the role of connecting link between the $\left( a\left( x \right),x \right)$, $\left\langle a\left( x \right),x \right\rangle $-algebras:
$$[\uparrow ,n]\left\langle x|a\left( x \right) \right\rangle ={{a}^{\left( n \right)}}\left( x \right), \qquad{{a}_{1}}=0.$$
For example,
$$\left\langle x|\frac{{{x}^{2}}}{1-x} \right\rangle =\left( \begin{matrix}
   0 & 0 & 0 & 0 & \ldots   \\
   1 & 0 & 0 & 0 & \ldots   \\
   0 & 1 & 0 & 0 & \ldots   \\
   0 & 1 & 0 & 0 & \ldots   \\
   0 & 1 & 1 & 0 & \ldots   \\
   0 & 1 & 0 & 0 & \ldots   \\
   0 & 1 & 2 & 0 & \ldots   \\
   0 & 1 & 0 & 0 & \ldots   \\
   0 & 1 & 2 & 1 & \ldots   \\
   0 & 1 & 1 & 0 & \ldots   \\
   0 & 1 & 2 & 0 & \ldots   \\
   0 & 1 & 0 & 0 & \ldots   \\
   0 & 1 & 4 & 3 & \ldots   \\
   \vdots  & \vdots  & \vdots  & \vdots  & \ddots   \\
\end{matrix} \right).$$
Denote
$$\left\langle x|a\left( x \right) \right\rangle f\left( x \right)=f\circ \left( a\left( x \right) \right),  \qquad f\left( x \right)=\sum\limits_{n=0}^{\infty }{{{f}_{n}}{{x}^{n}}};$$
$$\left\langle b\left( x \right),x \right\rangle \left\langle x|a\left( x \right) \right\rangle =\left\langle b\left( x \right)|a\left( x \right) \right\rangle .$$
Product of matrices of the form $\left\langle b\left( x \right)|a\left( x \right) \right\rangle $ is not a matrix of the same form, but since
$$\left\langle x|a\left( x \right) \right\rangle {{x}^{m}}f\left( x \right)={{a}^{\left( m \right)}}\left( x \right)\circ f\circ \left( a\left( x \right) \right),$$
then
$$\left\langle x|a\left( x \right) \right\rangle \left( f\left( x \right),x \right)=\left\langle f\circ \left( a\left( x \right) \right)|a\left( x \right) \right\rangle ,$$
$$\left\langle x|a\left( x \right) \right\rangle f\left( x \right)c\left( x \right)=f\circ \left( a\left( x \right) \right)\circ c\circ \left( a\left( x \right) \right),$$
$$\left\langle x|a\left( x \right) \right\rangle \left( 1,g\left( x \right) \right)=\left\langle x|g\circ \left( a\left( x \right) \right) \right\rangle ,$$
	$$\left\langle b\left( x \right)|a\left( x \right) \right\rangle \left( f\left( x \right),g\left( x \right) \right)=\left\langle b\left( x \right)\circ f\circ \left( a\left( x \right) \right)|g\circ \left( a\left( x \right) \right) \right\rangle. \eqno(1)$$

Thus, any matrix of the form $\left\langle b\left( x \right)|a\left( x \right) \right\rangle $ can be represented as the product of matrix of the same form and Riordan array. 

Here we get the definitions of the power and of the logarithm for the $\left\langle a\left( x \right),x \right\rangle $-algebra. Denote
$${{a}^{\left( \varphi  \right)}}\left( x \right)=\left\langle x|a\left( x \right)-x \right\rangle {{\left( 1+x \right)}^{\varphi }},  \qquad{{a}_{1}}=1.$$
Then
$${{a}^{\left( \varphi  \right)}}\left( x \right)\circ {{a}^{\left( \beta  \right)}}\left( x \right)={{a}^{\left( \varphi +\beta  \right)}}\left( x \right),\qquad{{\left( {{a}^{\left( \varphi  \right)}}\left( x \right) \right)}^{\left( \beta  \right)}}={{a}^{\left( \varphi \beta  \right)}}\left( x \right).$$
Denote
$$\log \circ a\left( x \right)=\left\langle x|a\left( x \right)-x \right\rangle \log \left( 1+x \right),  \qquad{{a}_{1}}=1.$$
Then
$$\log \circ {{a}^{\left( \varphi  \right)}}\left( x \right)=\varphi \log \circ a\left( x \right),\qquad\log \circ \left( a\left( x \right)\circ b\left( x \right) \right)=\log \circ a\left( x \right)+\log \circ b\left( x \right),$$
$$\left\langle x|\log \circ a\left( x \right) \right\rangle {{e}^{\varphi x}}={{a}^{\left( \varphi  \right)}}\left( x \right),\qquad{{a}^{\left( \varphi  \right)}}\left( x \right)\circ {{b}^{\left( \varphi  \right)}}\left( x \right)={{\left( a\left( x \right)\circ b\left( x \right) \right)}^{\left( \varphi  \right)}}.$$
Denote
$${{\left| {{e}^{x}} \right|}^{-1}}\left\langle x|\log \circ a\left( x \right) \right\rangle \left| {{e}^{x}} \right|={{\left\langle x|\log \circ a\left( x \right) \right\rangle }_{{{e}^{x}}}},\qquad\left[ n,\to  \right]{{\left\langle x|\log \circ a\left( x \right) \right\rangle }_{{{e}^{x}}}}={{\tilde{s}}_{n}}\left( x \right).$$
Then
$${{a}^{\left( \varphi  \right)}}\left( x \right)=\sum\limits_{n=0}^{\infty }{\frac{{{{\tilde{s}}}_{n}}\left( \varphi  \right)}{n!}}{{x}^{n}}.$$
\section{ }

We note the following analogy between the matrices $\left( 1,a\left( x \right) \right)$ and $\left\langle x|a\left( x \right) \right\rangle $. Denote

$${{B}_{n,m}}\left( {{a}_{1}},{{a}_{2}},...,{{a}_{n}} \right)=\sum{\frac{m!}{{{m}_{1}}!{{m}_{2}}!\text{ }...\text{ }{{m}_{n}}!}}\text{ }a_{1}^{{{m}_{1}}}a_{2}^{{{m}_{2}}}...\text{ }a_{n}^{{{m}_{n}}},   \qquad n>0,$$
where expression $\prod\nolimits_{k=1}^{n}{a_{k}^{{{m}_{k}}}}$ corresponding to the partition $n=\sum\nolimits_{k=1}^{n}{k{{m}_{k}}}$, $\sum\nolimits_{k=1}^{n}{{{m}_{k}}}=m$ and summation is done over all partitions of number $n$ into $m$ parts. Since
$${{\left( \sum\limits_{k=1}^{n}{{{a}_{k}}{{x}^{k}}} \right)}^{m}}=\sum\limits_{{{m}_{1}}+{{m}_{2}}+...+{{m}_{n}}=m}{\frac{m!}{{{m}_{1}}!{{m}_{2}}!...{{m}_{n}}!}}\text{ }{{\left( {{a}_{1}}x \right)}^{{{m}_{1}}}}{{\left( {{a}_{2}}{{x}^{2}} \right)}^{{{m}_{2}}}}...\text{ }{{\left( {{a}_{n}}{{x}^{n}} \right)}^{{{m}_{n}}}}=$$
$$=\sum\limits_{{{m}_{1}}+{{m}_{2}}+...+{{m}_{n}}=m}{\frac{m!}{{{m}_{1}}!{{m}_{2}}!...{{m}_{n}}!}}\text{ }a_{1}^{{{m}_{1}}}a_{2}^{{{m}_{2}}}...\text{ }a_{n}^{{{m}_{n}}}{{x}^{\sum }},   \qquad\sum ={{m}_{1}}+2{{m}_{2}}+...+n{{m}_{n}},$$
then
$$\left[ {{x}^{n}} \right]{{\left( \sum\limits_{n=1}^{\infty }{{{a}_{n}}{{x}^{n}}} \right)}^{m}}={{B}_{n,m}}\left( {{a}_{1}},{{a}_{2}},...,{{a}_{n}} \right),   \quad\left[ n,\to  \right]\left( 1,a\left( x \right) \right)=\sum\limits_{m=1}^{n}{{{B}_{n,m}}\left( {{a}_{1}},{{a}_{2}},...,{{a}_{n}} \right){{x}^{m}}}:$$
$$\left( 1,a\left( x \right) \right)=$$
$$\left( \begin{matrix}
   1 & 0 & 0 & 0 & 0 & 0 & \ldots   \\
   0 & {{a}_{1}} & 0 & 0 & 0 & 0 & \dots   \\
   0 & {{a}_{2}} & a_{1}^{2} & 0 & 0 & 0 & \dots   \\
   0 & {{a}_{3}} & 2{{a}_{1}}{{a}_{2}} & a_{1}^{3} & 0 & 0 &  \ldots   \\
   0 & {{a}_{4}} & 2{{a}_{1}}{{a}_{3}}+a_{2}^{2} & 3a_{1}^{2}{{a}_{2}} & a_{1}^{4} & 0 &  \ldots   \\
   0 & {{a}_{5}} & 2{{a}_{1}}{{a}_{4}}+2{{a}_{2}}{{a}_{3}} & 3a_{1}^{2}{{a}_{3}}+3{{a}_{1}}a_{2}^{2} & 4a_{1}^{3}{{a}_{2}} & a_{1}^{5} & \ldots   \\
   0 & {{a}_{6}} & 2{{a}_{1}}{{a}_{5}}+2{{a}_{2}}{{a}_{4}}+a_{3}^{2} & 3a_{1}^{2}{{a}_{4}}+6{{a}_{1}}{{a}_{2}}{{a}_{3}}+a_{2}^{3} & 4a_{1}^{3}{{a}_{3}}+6a_{1}^{2}a_{2}^{2} & 5a_{1}^{4}{{a}_{2}} & \ldots   \\
   \vdots  & \vdots  & \vdots  & \vdots  & \vdots  & \vdots  &  \ddots   \\
\end{matrix} \right)$$

If ${{a}_{0}}=1$, $\log a\left( x \right)=b\left( x \right)$, ${{s}_{n}}\left( x \right)=\left[ n,\to  \right]{{\left( 1,b\left( x \right) \right)}_{{{e}^{x}}}}$, then 
$${{s}_{n}}\left( x \right)=n!\sum\limits_{m=1}^{n}{\frac{{{B}_{n,m}}\left( {{b}_{1}},{{b}_{2}},...,{{b}_{n}} \right)}{m!}{{x}^{m}}},  \qquad{{b}_{n}}=\sum\limits_{m=1}^{n}{{{\left( -1 \right)}^{m+1}}\frac{{{B}_{n,m}}\left( {{a}_{1}},{{a}_{2}},...,{{a}_{n}} \right)}{m}}.$$

Denote
$${{\tilde{B}}_{n,m}}\left( {{a}_{2}},{{a}_{3}},...,{{a}_{n}} \right)=\sum{\frac{m!}{{{m}_{2}}!{{m}_{3}}!\text{ }...\text{ }{{m}_{n}}!}}\text{ }a_{2}^{{{m}_{2}}}a_{3}^{{{m}_{3}}}...\text{ }a_{n}^{{{m}_{n}}},  \qquad n>1,$$
where expression $\prod\nolimits_{k=2}^{n}{a_{k}^{{{m}_{k}}}}$ corresponding to the decomposition $n=\prod\nolimits_{k=2}^{n}{{{k}^{{{m}_{k}}}}}$, ${\sum\nolimits_{k=2}^{n}{{{m}_{k}}}=m}$, and summation is done over all decompositions of number $n$ into $m$ factors . Since
$${{\left( \sum\limits_{k=2}^{n}{{{a}_{k}}{{x}^{k}}} \right)}^{\left( m \right)}}=$$
$$=\sum\limits_{{{m}_{2}}+{{m}_{3}}+...+{{m}_{n}}=m}{\frac{m!}{{{m}_{2}}!{{m}_{3}}!...{{m}_{n}}!}}\text{ }{{\left( {{a}_{2}}{{x}^{2}} \right)}^{\left( {{m}_{2}} \right)}}\circ {{\left( {{a}_{3}}{{x}^{3}} \right)}^{\left( {{m}_{3}} \right)}}\circ ...\circ {{\left( {{a}_{n}}{{x}^{n}} \right)}^{\left( {{m}_{n}} \right)}}=$$
$$=\sum\limits_{{{m}_{2}}+{{m}_{3}}+...+{{m}_{n}}=m}{\frac{m!}{{{m}_{2}}!{{m}_{3}}!...{{m}_{n}}!}}\text{ }a_{2}^{{{m}_{2}}}a_{3}^{{{m}_{3}}}...\text{ }a_{n}^{{{m}_{n}}}{{x}^{\prod }},  \qquad\prod ={{2}^{{{m}_{2}}}}{{3}^{{{m}_{3}}}}...\text{ }{{n}^{{{m}_{n}}}},$$
then
$$\left[ {{x}^{n}} \right]{{\left( \sum\limits_{n=2}^{\infty }{{{a}_{n}}{{x}^{n}}} \right)}^{\left( m \right)}}={{\tilde{B}}_{n,m}}\left( {{a}_{2}},{{a}_{3}},...,{{a}_{n}} \right),   \quad\left[ n,\to  \right]\left\langle x|a\left( x \right) \right\rangle =\sum\limits_{m=1}^{n}{{{{\tilde{B}}}_{n,m}}\left( {{a}_{2}},{{a}_{3}},...,{{a}_{n}} \right){{x}^{m}}}:$$
$$\left\langle x|a\left( x \right) \right\rangle =\left( \begin{matrix}
   0 & 0 & 0 & 0 & 0 & \ldots   \\
   1 & 0 & 0 & 0 & 0 & \ldots   \\
   0 & {{a}_{2}} & 0 & 0 & 0 & \ldots   \\
   0 & {{a}_{3}} & 0 & 0 & 0 & \ldots   \\
   0 & {{a}_{4}} & a_{2}^{2} & 0 & 0 & \ldots   \\
   0 & {{a}_{5}} & 0 & 0 & 0 & \ldots   \\
   0 & {{a}_{6}} & 2{{a}_{2}}{{a}_{3}} & 0 & 0 & \ldots   \\
   0 & {{a}_{7}} & 0 & 0 & 0 & \ldots   \\
   0 & {{a}_{8}} & 2{{a}_{2}}{{a}_{4}} & a_{2}^{3} & 0 & \ldots   \\
   0 & {{a}_{9}} & a_{3}^{2} & 0 & 0 & \ldots   \\
   0 & {{a}_{10}} & 2{{a}_{2}}{{a}_{5}} & 0 & 0 & \ldots   \\
   0 & {{a}_{11}} & 0 & 0 & 0 & \ldots   \\
   0 & {{a}_{12}} & 2{{a}_{2}}{{a}_{6}}+2{{a}_{4}}{{a}_{3}} & 3a_{2}^{2}{{a}_{3}} & 0 & \ldots   \\
   0 & {{a}_{13}} & 0 & 0 & 0 & \ldots   \\
   0 & {{a}_{14}} & 2{{a}_{2}}{{a}_{7}} & 0 & 0 & \ldots   \\
   0 & {{a}_{15}} & 2{{a}_{3}}{{a}_{5}} & 0 & 0 & \ldots   \\
   0 & {{a}_{16}} & 2{{a}_{2}}{{a}_{8}}+a_{4}^{2} & 3a_{2}^{2}{{a}_{4}} & a_{2}^{4} & \ldots   \\
   \vdots  & \vdots  & \vdots  & \vdots  & \vdots  & \ddots   \\
\end{matrix} \right)$$

If ${{a}_{1}}=1$, $\log \circ a\left( x \right)=b\left( x \right)$, ${{\tilde{s}}_{n}}\left( x \right)=\left[ n,\to  \right]{{\left\langle x|b\left( x \right) \right\rangle }_{{{e}^{x}}}}$, then 
$${{\tilde{s}}_{n}}\left( x \right)=n!\sum\limits_{m=1}^{n}{\frac{{{{\tilde{B}}}_{n,m}}\left( {{b}_{2}},{{b}_{3}},...,{{b}_{n}} \right)}{m!}{{x}^{m}}},    \qquad{{b}_{n}}=\sum\limits_{m=1}^{n}{{{\left( -1 \right)}^{m+1}}\frac{{{{\tilde{B}}}_{n,m}}\left( {{a}_{2}},{{a}_{3}},...,{{a}_{n}} \right)}{m}}.$$

\section{ }

Relationship between the $\left( a\left( x \right),x \right)$, $\left\langle a\left( x \right),x \right\rangle $-algebras is most visibly manifested in the following theorem.\\
{\bfseries Theorem 1.}\emph{ Each formal power series $a\left( x \right)$, ${{a}_{0}}=1$, 
$${{a}^{\varphi }}\left( x \right)=\sum\limits_{n=0}^{\infty }{\frac{{{s}_{n}}\left( \varphi  \right)}{n!}{{x}^{n}}}, \qquad{{s}_{n}}\left( x \right)=\left[ n,\to  \right]{{\left( 1,\log a\left( x \right) \right)}_{{{e}^{x}}}},$$
corresponds to the series $\tilde{a}\left( x \right)$,
$${{\tilde{a}}^{\left( \varphi  \right)}}\left( x \right)=x+\sum\limits_{n=2}^{\infty }{\frac{{{s}_{{{m}_{1}}}}\left( \varphi  \right){{s}_{{{m}_{2}}}}\left( \varphi  \right)...{{s}_{{{m}_{r}}}}\left( \varphi  \right)}{{{m}_{1}}!{{m}_{2}}!\text{ }...\text{ }{{m}_{r}}!}}{{x}^{n}},  \qquad n=p_{1}^{{{m}_{1}}}p_{2}^{{{m}_{2}}}...\text{ }p_{r}^{{{m}_{r}}},$$
where $n=p_{1}^{{{m}_{1}}}p_{2}^{{{m}_{2}}}...\text{ }p_{r}^{{{m}_{r}}}$ is the canonical decomposition of number $n$.}\\
{\bfseries Proof.} We will denote a prime number of the letter $p$. Consider series
$${{\tilde{a}}^{\left( \varphi  \right)}}\left( x \right)=\prod\limits_{p=2}^{\infty }{\circ {{\left( \sum\limits_{n=0}^{\infty }{{{a}_{p,n}}{{x}^{{{p}^{n}}}}} \right)}^{\left( \varphi  \right)}}}=\prod\limits_{p=2}^{\infty }{\circ \tilde{a}_{p}^{\left( \varphi  \right)}}\left( x \right), \qquad{{a}_{p,0}}=1,$$
where product, denoted similar to the ordinary product, is taken over all prime numbers. In view of the isomorphism between the algebra of matrices $\left\langle a\left( x \right),x \right\rangle =\sum\nolimits_{n=0}^{\infty }{{{a}_{n}}}{{\left\langle {{x}^{p}},x \right\rangle }^{n}}$ and the algebra of matrices $\left( a\left( x \right),x \right)$, the series ${{\tilde{a}}_{p}}\left( x \right)$ corresponds to the series ${{a}_{p}}\left( x \right)$, such that
$$\left[ {{x}^{{{p}^{n}}}} \right]\tilde{a}_{p}^{\left( \varphi  \right)}\left( x \right)=\left[ {{x}^{n}} \right]a_{p}^{\varphi }\left( x \right)=\frac{{{s}_{p,n}}\left( \varphi  \right)}{n!}, \qquad{{s}_{p,n}}\left( x \right)=\left[ n,\to  \right]{{\left( 1,\log {{a}_{p}}\left( x \right) \right)}_{{{e}^{x}}}}.$$
Hence,
$${{\tilde{a}}^{\left( \varphi  \right)}}\left( x \right)=x+\sum\limits_{n=2}^{\infty }{\frac{{{s}_{{{p}_{1}},{{m}_{1}}}}\left( \varphi  \right){{s}_{{{p}_{2}},{{m}_{2}}}}\left( \varphi  \right)...{{s}_{{{p}_{_{r}}},{{m}_{r}}}}\left( \varphi  \right)}{{{m}_{1}}!{{m}_{2}}!\text{ }...\text{ }{{m}_{r}}!}}{{x}^{n}},  \qquad n=p_{1}^{{{m}_{1}}}p_{2}^{{{m}_{2}}}...\text{ }p_{r}^{{{m}_{r}}}.$$
We are interested in the case when all the series ${{\tilde{a}}_{p}}\left( x \right)$ corresponds to the same series $a\left( x \right)$. In this case the series $\tilde{a}\left( x \right)$ form a group, isomorphic to a group of  the series $a\left( x \right)$: if $a\left( x \right)b\left( x \right)=c\left( x \right)$, then $\tilde{a}\left( x \right)\circ \tilde{b}\left( x \right)=\tilde{c}\left( x \right)$. Note that
$$\left[ {{x}^{n}} \right]\log \circ \tilde{a}\left( x \right)=0, \qquad n\ne {{p}^{m}};  \qquad=\left[ {{x}^{m}} \right]\log a\left( x \right), \qquad n={{p}^{m}}.$$

Denote $\zeta \left( x \right)=\sum\nolimits_{n=1}^{\infty }{{{x}^{n}}}$. Then
$${{\varsigma }^{\left( \varphi  \right)}}\left( x \right)=x+\sum\limits_{n=2}^{\infty }{\frac{{{\left( \varphi  \right)}^{{{m}_{1}}}}{{\left( \varphi  \right)}^{{{m}_{2}}}}...{{\left( \varphi  \right)}^{{{m}_{r}}}}}{{{m}_{1}}!{{m}_{2}}!\text{ }...\text{ }{{m}_{r}}!}{{x}^{n}}}, \qquad n=p_{1}^{{{m}_{1}}}p_{2}^{{{m}_{2}}}...\text{ }p_{r}^{{{m}_{r}}},$$
where ${{\left( \varphi  \right)}^{{{m}_{i}}}}=\varphi \left( \varphi +1 \right)\left( \varphi +2 \right)...\left( \varphi +{{m}_{i}}-1 \right)$. Denote
$$\left[ n,\to  \right]{{\left\langle x|\log \circ \varsigma \left( x \right) \right\rangle }_{{{e}^{x}}}}={{\tilde{s}}_{n}}\left( x \right).$$
Then
$${{\tilde{s}}_{0}}\left( x \right)=0,   \qquad{{\tilde{s}}_{1}}\left( x \right)=1,    \qquad\frac{{{{\tilde{s}}}_{n}}\left( x \right)}{n!}=\frac{{{\left( x \right)}^{{{m}_{1}}}}{{\left( x \right)}^{{{m}_{2}}}}...{{\left( x \right)}^{{{m}_{r}}}}}{{{m}_{1}}!{{m}_{2}}!\text{ }...\text{ }{{m}_{r}}!}.$$
A-priory,
$$\left[ {{x}^{n}} \right]\log \circ \varsigma \left( x \right)=0, \qquad n\ne {{p}^{m}};  \qquad=\frac{1}{m}, \qquad n={{p}^{m}}.$$
On the other hand,
$$\left[ {{x}^{n}} \right]\log \circ \varsigma \left( x \right)=\sum\limits_{m=1}^{n}{{{\left( -1 \right)}^{m+1}}\frac{{{{\tilde{B}}}_{n,m}}\left( 1,1,...,1 \right)}{m}},$$ 
where 
$${{\tilde{B}}_{n,m}}\left( 1,1,...,1 \right)=\sum{\frac{m!}{{{m}_{2}}!{{m}_{3}}!\text{ }...\text{ }{{m}_{n}}!}},    \qquad n=\prod\limits_{k=2}^{n}{{{k}^{{{m}_{k}}}}},  \qquad\sum\limits_{k=2}^{n}{{{m}_{k}}}=m,$$
and summation is done over all decompositions of number $n$ into $m$ factors. We note also the identity 
$$\sum\limits_{m=1}^{n}{\left( \begin{matrix}
   \varphi   \\
   m  \\
\end{matrix} \right){{{\tilde{B}}}_{n,m}}}\left( 1,1,...,1 \right)=\left( \begin{matrix}
   \varphi +{{s}_{1}}-1  \\
   {{s}_{1}}  \\
\end{matrix} \right)\left( \begin{matrix}
   \varphi +{{s}_{2}}-1  \\
   {{s}_{2}}  \\
\end{matrix} \right)...\left( \begin{matrix}
   \varphi +{{s}_{r}}-1  \\
   {{s}_{r}}  \\
\end{matrix} \right),$$
$n=p_{1}^{{{s}_{1}}}p_{2}^{{{s}_{2}}}...\text{ }p_{r}^{{{s}_{r}}}$, similar to the identity
$$\sum\limits_{m=1}^{n}{\left( \begin{matrix}
   \varphi   \\
   m  \\
\end{matrix} \right){{B}_{n,m}}}\left( 1,1,...,1 \right)=\left( \begin{matrix}
   \varphi +n-1  \\
   n  \\
\end{matrix} \right),$$
where
$${{B}_{n,m}}\left( 1,1,...,1 \right)=\left( \begin{matrix}
   n-1  \\
   m-1  \\
\end{matrix} \right)=\sum{\frac{m!}{{{m}_{1}}!{{m}_{2}}!\text{ }...\text{ }{{m}_{n}}!}}\text{ }, \qquad n=\sum\limits_{k=1}^{n}{k{{m}_{k}}}, \qquad\sum\limits_{k=1}^{n}{{{m}_{k}}}=m,$$
and summation is done over all partitions of number $n$ into $m$ parts. 

Analog of the exponential series for the $\left\langle a\left( x \right),x \right\rangle $-algebra is the series
$${{\varepsilon }^{\left( \varphi  \right)}}\left( x \right)=x+\sum\limits_{n=2}^{\infty }{\frac{{{\varphi }^{{{m}_{1}}+{{m}_{2}}+...+{{m}_{r}}}}}{{{m}_{1}}!{{m}_{2}}!\text{ }...\text{ }{{m}_{r}}!}{{x}^{n}}},  \qquad n=p_{1}^{{{m}_{1}}}p_{2}^{{{m}_{2}}}...\text{ }p_{r}^{{{m}_{r}}}.$$
 Series $\log \circ \varepsilon \left( x \right)$ is closely connected with the sequence of prime numbers:
$$\left[ {{x}^{n}} \right]\log \circ \varepsilon \left( x \right)=0,  \qquad n\ne p; \qquad=1,  \qquad n=p.$$
In general  case 
$${{\left( \log \circ \varepsilon \left( x \right) \right)}^{\left( m \right)}}=\sum{\frac{m!}{{{m}_{1}}!{{m}_{2}}!...{{m}_{r}}!}}{{x}^{n}},$$
where  summation is done over all $n=p_{1}^{{{m}_{1}}}p_{2}^{{{m}_{2}}}...\text{ }p_{r}^{{{m}_{r}}}$, ${{m}_{1}}+{{m}_{2}}+...+{{m}_{r}}=m$. 

\section{ }

It follows from the Lagrange series expansion  for arbitrary formal power series $b\left( x \right)$ and $a\left( x \right)$, ${{a}_{0}}=1$:
$$\frac{b\left( x \right)}{1-x{{\left( \log a\left( x \right) \right)}^{\prime }}}=\sum\limits_{n=0}^{\infty }{\frac{{{x}^{n}}}{{{a}^{n}}\left( x \right)}}\left[ {{x}^{n}} \right]b\left( x \right){{a}^{n}}\left( x \right)$$
that each formal power series $a\left( x \right)$, ${{a}_{0}}=1$, is associated by means of the transform
$${{a}^{\varphi }}\left( x \right)=\sum\limits_{n=0}^{\infty }{\frac{{{x}^{n}}}{{{a}^{\beta n}}\left( x \right)}\left[ {{x}^{n}} \right]}\left( 1-x\beta {{\left( \log a\left( x \right) \right)}^{\prime }} \right){{a}^{\varphi +\beta n}}\left( x \right)$$
with the set of series $_{\left( \beta  \right)}a\left( x \right)$, $_{\left( 0 \right)}a\left( x \right)=a\left( x \right)$,   such that
$${}_{\left( \beta  \right)}a\left( x{{a}^{-\beta }}\left( x \right) \right)=a\left( x \right),   \qquad a\left( x{}_{\left( \beta  \right)}{{a}^{\beta }}\left( x \right) \right)={}_{\left( \beta  \right)}a\left( x \right),$$
$$\left[ {{x}^{n}} \right]{}_{\left( \beta  \right)}{{a}^{\varphi }}\left( x \right)=\left[ {{x}^{n}} \right]\left( 1-x\beta \frac{{a}'\left( x \right)}{a\left( x \right)} \right){{a}^{\varphi +\beta n}}\left( x \right)=\frac{\varphi }{\varphi +\beta n}\left[ {{x}^{n}} \right]{{a}^{\varphi +\beta n}}\left( x \right),$$
$$\left[ {{x}^{n}} \right]\left( 1+x\beta \frac{_{\left( \beta  \right)}{a}'\left( x \right)}{_{\left( \beta  \right)}a\left( x \right)} \right){}_{\left( \beta  \right)}{{a}^{\varphi }}\left( x \right)=\frac{\varphi +\beta n}{\varphi }\left[ {{x}^{n}} \right]{}_{\left( \beta  \right)}{{a}^{\varphi }}\left( x \right)=\left[ {{x}^{n}} \right]a_{{}}^{\varphi +\beta n}\left( x \right).$$
$${{\left( 1,x{}_{\left( \beta  \right)}{{a}^{\beta }}\left( x \right) \right)}^{-1}}=\left( 1,x{{a}^{-\beta }}\left( x \right) \right),$$
$${{\left( 1+x\beta {{\left( \log {}_{\left( \beta  \right)}a\left( x \right) \right)}^{\prime }},x{}_{\left( \beta  \right)}{{a}^{\beta }}\left( x \right) \right)}^{-1}}=\left( 1-x\beta {{\left( \log a\left( x \right) \right)}^{\prime }},x{{a}^{-\beta }}\left( x \right) \right),$$
$$\left[ n,\to  \right]\left( 1,x{}_{\left( \beta  \right)}{{a}^{\beta }}\left( x \right) \right)=\left[ n,\to  \right]\left( 1-x\beta {{\left( \log a\left( x \right) \right)}^{\prime }}{{a}^{\beta n}}\left( x \right),x \right),$$
$$\left[ n,\to  \right]\left( 1+x\beta {{\left( \log {}_{\left( \beta  \right)}a\left( x \right) \right)}^{\prime }},x{}_{\left( \beta  \right)}{{a}^{\beta }}\left( x \right) \right)=\left[ n,\to  \right]\left( {{a}^{\beta n}}\left( x \right),x \right).$$
Denote
$$\left[ n,\to  \right]{{\left( 1,\log {}_{\left( \beta  \right)}a\left( x \right) \right)}_{{{e}^{x}}}}={}_{\left( \beta  \right)}{{s}_{n}}\left( x \right),  \qquad_{\left( 0 \right)}{{s}_{n}}\left( x \right)={{s}_{n}}\left( x \right).$$
Then
$$_{\left( \beta  \right)}{{a}^{\varphi }}\left( x \right)=\sum\limits_{n=0}^{\infty }{\frac{\varphi }{\varphi +\beta n}}\frac{{{s}_{n}}\left( \varphi +\beta n \right)}{n!}{{x}^{n}},\qquad_{\left( \beta  \right)}{{s}_{n}}\left( x \right)=x{{\left( x+\beta n \right)}^{-1}}{{s}_{n}}\left( x+\beta n \right).$$

Apparently, the series $_{\left( \beta  \right)}a\left( x \right)$ for integer $\beta $, denoted by ${{S}_{\beta }}\left( x \right)$, were first considered in [6]. In [7] these series, called generalized Lagrange series, are considered in connection with the Riordan arrays. Examples of this construction are the generalized binomial and generalized exponential series [8; p. 200]:
$$a\left( x \right)=1+x,\qquad{}_{\left( \beta  \right)}{{a}^{\varphi }}\left( x \right)=\sum\limits_{n=0}^{\infty }{\frac{\varphi }{\varphi +\beta n}}\left( \begin{matrix}
   \varphi +\beta n  \\
   n  \\
\end{matrix} \right){{x}^{n}};$$
$$a\left( x \right)={{e}^{x}},  \qquad{}_{\left( \beta  \right)}{{a}^{\varphi }}\left( x \right)=\sum\limits_{n=0}^{\infty }{\frac{\varphi {{\left( \varphi +\beta n \right)}^{n-1}}}{n!}{{x}^{n}}}.$$

We introduce analog of the differential operator for the $\left\langle a\left( x \right),x \right\rangle $-algebra:
$$\tilde{D}a\left( x \right)={{a}^{*}}\left( x \right)=\sum\limits_{n=1}^{\infty }{\ln n{{a}_{n}}{{x}^{n}}}.$$
Since
$$\ln n\sum\limits_{d|n}{{{a}_{d}}}{{b}_{{n}/{d}\;}}=\sum\limits_{d|n}{\ln \left( {n}/{d}\; \right){{a}_{d}}{{b}_{{n}/{d}\;}}}+\sum\limits_{d|n}{\ln d{{a}_{d}}{{b}_{{n}/{d}\;}}},$$
then
$${{\left( a\left( x \right)\circ b\left( x \right) \right)}^{*}}=a\left( x \right)\circ {{b}^{*}}\left( x \right)+{{a}^{*}}\left( x \right)\circ b\left( x \right), \qquad{{\left( {{a}^{\left( n \right)}}\left( x \right) \right)}^{*}}=n{{a}^{\left( n-1 \right)}}\left( x \right)\circ {{a}^{*}}\left( x \right),$$
$${{\left( \log \circ a\left( x \right) \right)}^{*}}={{a}^{*}}\left( x \right)\circ \sum\limits_{n=1}^{\infty }{{{\left( -1 \right)}^{n-1}}{{\left( a\left( x \right)-x \right)}^{\left( n-1 \right)}}}={{a}^{*}}\left( x \right)\circ {{a}^{\left( -1 \right)}}\left( x \right).$$
Since
$${{a}^{\left( \varphi  \right)}}\left( x \right)=\left\langle x|a\left( x \right)-x \right\rangle {{\left( 1+x \right)}^{\varphi }}, \qquad\tilde{D}\left\langle x|a\left( x \right)-x \right\rangle =\left\langle {{a}^{*}}\left( x \right)|a\left( x \right)-x \right\rangle D,$$
where $D$ is the matrix of  differential operator, then
$${{\left( {{a}^{\left( \varphi  \right)}}\left( x \right) \right)}^{*}}=\varphi {{a}^{\left( \varphi -1 \right)}}\left( x \right)\circ {{a}^{*}}\left( x \right).$$
{\bfseries Theorem 2.}\emph { Each formal power series $a\left( x \right)$, ${{a}_{0}}=0$, ${{a}_{1}}=1$, is associated by means of the transform
$${{a}^{\left( \varphi  \right)}}\left( x \right)=\sum\limits_{n=1}^{\infty }{{{a}^{\left( -\beta \ln n \right)}}\left( {{x}^{n}} \right)\left[ {{x}^{n}} \right]}\left( x-\beta {{\left( \log \circ a\left( x \right) \right)}^{*}} \right)\circ {{a}^{\left( \varphi +\beta \ln n \right)}}\left( x \right)$$
with the set of series $_{\left( \beta  \right)}a\left( x \right)$, $_{\left( 0 \right)}a\left( x \right)=a\left( x \right)$,   such that}
$$\left[ {{x}^{n}} \right]{}_{\left( \beta  \right)}{{a}^{\left( \varphi  \right)}}\left( x \right)=\left[ {{x}^{n}} \right]\left( x-\beta {{a}^{*}}\left( x \right)\circ {{a}^{\left( -1 \right)}}\left( x \right) \right)\circ {{a}^{\left( \varphi +\beta \ln n \right)}}\left( x \right)=$$
$$=\frac{\varphi }{\varphi +\beta \ln n}\left[ {{x}^{n}} \right]{{a}^{\left( \varphi +\beta \ln n \right)}}\left( x \right),$$
$$\left[ {{x}^{n}} \right]\left( x+\beta {}_{\left( \beta  \right)}{{a}^{*}}\left( x \right)\circ {}_{\left( \beta  \right)}{{a}^{\left( -1 \right)}}\left( x \right) \right)\circ {}_{\left( \beta  \right)}{{a}^{\left( \varphi  \right)}}\left( x \right)=\frac{\varphi +\beta \ln n}{\varphi }\left[ {{x}^{n}} \right]{}_{\left( \beta  \right)}{{a}^{\left( \varphi  \right)}}\left( x \right)=$$
$$=\left[ {{x}^{n}} \right]{{a}^{\left( \varphi +\beta \ln n \right)}}\left( x \right).$$

For proof we introduce the matrices $\left\langle x,a\left( x \right) \right\rangle $:
$$\left\langle x,a\left( x \right) \right\rangle =\left( \begin{matrix}
   {{\left( {{a}_{\Sigma }} \right)}^{\ln 0}} & 0 & 0 & 0 & 0 & 0 & 0 & 0 & 0 & 0 & \ldots   \\
   0 & 1 & 0 & 0 & 0 & 0 & 0 & 0 & 0 & 0 & \ldots   \\
   0 & 0 & a_{1}^{2} & 0 & 0 & 0 & 0 & 0 & 0 & 0 & \ldots   \\
   0 & 0 & 0 & a_{1}^{3} & 0 & 0 & 0 & 0 & 0 & 0 & \ldots   \\
   0 & 0 & a_{2}^{2} & 0 & a_{1}^{4} & 0 & 0 & 0 & 0 & 0 & \ldots   \\
   0 & 0 & 0 & 0 & 0 & a_{1}^{5} & 0 & 0 & 0 & 0 & \ldots   \\
   0 & 0 & a_{3}^{2} & a_{2}^{3} & 0 & 0 & a_{1}^{6} & 0 & 0 & 0 & \ldots   \\
   0 & 0 & 0 & 0 & 0 & 0 & 0 & a_{1}^{7} & 0 & 0 & \ldots   \\
   0 & 0 & a_{4}^{2} & 0 & a_{2}^{4} & 0 & 0 & 0 & a_{1}^{8} & 0 & \ldots   \\
   0 & 0 & 0 & a_{3}^{3} & 0 & 0 & 0 & 0 & 0 & a_{1}^{9} & \ldots   \\
   \vdots  & \vdots  & \vdots  & \vdots  & \vdots  & \vdots  & \vdots  & \vdots  & \vdots  & \vdots  & \ddots   \\
\end{matrix} \right),$$ 
$$a_{n}^{m}=\left[ {{x}^{n}} \right]{{a}^{\left( \ln m \right)}}\left( x \right),   \qquad[\uparrow ,n]\left\langle x,a\left( x \right) \right\rangle ={{a}^{\left( \ln n \right)}}\left( {{x}^{n}} \right).$$
Denote
$$\left\langle x,a\left( x \right) \right\rangle b\left( x \right)={{b}_{d}}\circ \left( a\left( x \right) \right),  \qquad\left[ {{x}^{n}} \right]{{b}_{d}}\circ \left( a\left( x \right) \right)=\sum\limits_{d|n}{{{b}_{d}}a_{{n}/{d}\;}^{d}}; $$
$$\left\langle b\left( x \right),x \right\rangle \left\langle x,a\left( x \right) \right\rangle =\left\langle b\left( x \right),a\left( x \right) \right\rangle .$$
Since 
$${{a}^{\left( \ln n \right)}}\left( {{x}^{n}} \right)={{x}^{n}}\circ {{a}^{\left( \ln n \right)}}\left( x \right), \qquad n>0,$$
Then
$$[\uparrow ,n]\left\langle b\left( x \right),a\left( x \right) \right\rangle ={{x}^{n}}\circ b\left( x \right)\circ {{a}^{\left( \ln n \right)}}\left( x \right).$$
Let ${{b}_{d}}\circ \left( a\left( x \right) \right)=c\left( x \right)$. If we accept the rules
$${{\left( {{a}_{\Sigma }} \right)}^{\ln 0}}{{b}_{\Sigma }}={{c}_{\Sigma }}{{\left( {{a}_{\Sigma }} \right)}^{\ln 0}}, \qquad{{\left( {{a}_{\Sigma }} \right)}^{\ln 0}}{{\left( {{b}_{\Sigma }} \right)}^{\ln 0}}={{\left( {{a}_{\Sigma }}{{c}_{\Sigma }} \right)}^{\ln 0}}, \qquad{{\left( 1 \right)}^{\ln 0}}=1,$$
then following theorem is true.\\
{\bfseries Theorem 3.}\emph{ Matrices $\left\langle b\left( x \right),a\left( x \right) \right\rangle $, ${{b}_{1}}\ne 0$, ${{a}_{1}}\ne 0$, form a group whose elements are multiplied by the rule}
$$\left\langle b\left( x \right),a\left( x \right) \right\rangle \left\langle f\left( x \right),g\left( x \right) \right\rangle =\left\langle b\left( x \right)\circ {{f}_{d}}\circ \left( a\left( x \right) \right),a\left( x \right)\circ {{g}_{d}}\circ \left( a\left( x \right) \right) \right\rangle .$$
{\bfseries Proof.} Since
$$\left\langle x,a\left( x \right) \right\rangle b\left( {{x}^{m}} \right)=\sum\limits_{n=1}^{\infty }{{{b}_{n}}{{a}^{\left( \ln mn \right)}}\left( {{x}^{mn}} \right)}={{x}^{m}}\circ {{a}^{\left( \ln m \right)}}\left( x \right)\circ \sum\limits_{n=1}^{\infty }{{{b}_{n}}}{{a}^{\left( \ln n \right)}}\left( {{x}^{n}} \right),$$
then 
$$\left\langle x,a\left( x \right) \right\rangle \left\langle b\left( x \right),x \right\rangle =\left\langle {{b}_{d}}\circ \left( a\left( x \right) \right),a\left( x \right) \right\rangle .$$
Thus,
$$\left\langle x,a\left( x \right) \right\rangle b\left( x \right)\circ c\left( x \right)={{b}_{d}}\circ \left( a\left( x \right) \right)\circ {{c}_{d}}\circ \left( a\left( x \right) \right).$$
Since
$$\left\langle x,a\left( x \right) \right\rangle {{b}^{\left( \ln m \right)}}\left( {{x}^{m}} \right)={{x}^{m}}\circ {{a}^{\left( \ln m \right)}}\left( x \right)\circ {{\left( {{b}_{d}}\circ \left( a\left( x \right) \right) \right)}^{\left( \ln m \right)}},$$
then
$$\left\langle x,a\left( x \right) \right\rangle \left\langle x,b\left( x \right) \right\rangle =\left\langle x,a\left( x \right)\circ {{b}_{d}}\circ \left( a\left( x \right) \right) \right\rangle .$$

As we shall see, with respect to the some structure that in the ordinary algebra of formal power series corresponds to the Lagrange series, a complete analogy exists between the group of matrices $\left\langle b\left( x \right),a\left( x \right) \right\rangle $, ${{b}_{1}}\ne 0$, ${{a}_{1}}\ne 0$, and the Riordan  group. So we call this group the Riordan-Dirichlet group.

Note the identity for the matrices $\left\langle b\left( x \right)|a\left( x \right) \right\rangle $, complementary to identity (1):
$$\left\langle f\left( x \right),g\left( x \right) \right\rangle \left\langle b\left( x \right)|a\left( x \right) \right\rangle =\left\langle f\left( x \right)\circ {{b}_{d}}\circ \left( g\left( x \right) \right)|{{a}_{d}}\circ \left( g\left( x \right) \right) \right\rangle .$$

Now we prove the theorem 2. Let the matrices $\left\langle x,{{a}^{\left( -1 \right)}}\left( x \right) \right\rangle $, $\left\langle x,b\left( x \right) \right\rangle $ are mutually inverse. Then
$$\left\langle x,{{a}^{\left( -1 \right)}}\left( x \right) \right\rangle b\left( x \right)=a\left( x \right),   \qquad\left\langle x,b\left( x \right) \right\rangle a\left( x \right)=b\left( x \right).$$
Let $\tilde{D}1=0$ (perhaps we should accept $\tilde{D}1=\ln 0$, but this is not fundamentally now). Since
$$\tilde{D}{{b}^{\left( \ln n \right)}}\left( {{x}^{n}} \right)={{\left( {{x}^{n}}\circ {{b}^{\left( \ln n \right)}}\left( x \right) \right)}^{*}}=
\ln n{{x}^{n}}\circ {{b}^{\left( \ln n \right)}}\left( x \right)\circ \left( x+{{b}^{*}}\left( x \right)\circ {{b}^{\left( -1 \right)}}\left( x \right) \right),$$
then
$$\tilde{D}\left\langle x,b\left( x \right) \right\rangle =\left\langle x+{{\left( \log \circ b\left( x \right) \right)}^{*}},b\left( x \right) \right\rangle \tilde{D},
\qquad\left\langle x+{{\left( \log \circ b\left( x \right) \right)}^{*}},b\left( x \right) \right\rangle {{a}^{*}}\left( x \right)={{b}^{*}}\left( x \right),$$
$${{\left\langle x+{{\left( \log \circ b\left( x \right) \right)}^{*}},b\left( x \right) \right\rangle }^{-1}}=\left\langle x-{{\left( \log \circ a\left( x \right) \right)}^{*}},{{a}^{\left( -1 \right)}}\left( x \right) \right\rangle .$$
Denote
$$\left[ {{x}^{n}} \right]{{a}^{\left( \ln m \right)}}\left( x \right)=a_{n}^{m},    \qquad\left[ {{x}^{n}} \right]\left( x-{{\left( \log \circ a\left( x \right) \right)}^{*}} \right)\circ {{a}^{\left( \ln m \right)}}\left( x \right)=c_{n}^{m},$$
$${{a}_{m}}\left( x \right)=\sum\limits_{n=1}^{\infty }{a_{n}^{mn}}{{x}^{n}},   \qquad{{c}_{m}}\left( x \right)=\sum\limits_{n=1}^{\infty }{c_{n}^{mn}}{{x}^{n}}.$$
 Construct the matrices $A$, $C$:
$$[\uparrow ,0]A=\varphi ,  \qquad[\uparrow ,n]A={{a}_{n}}\left( {{x}^{n}} \right),  \qquad[\uparrow ,0]C=\varphi ,  \qquad[\uparrow ,n]C={{c}_{n}}\left( {{x}^{n}} \right),$$
where $\varphi $ is a some number,
$$A=\left( \begin{matrix}
   \varphi& 0 & 0 & 0 & 0 & 0 & 0 & \ldots   \\
   0 & a_{1}^{1} & 0 & 0 & 0 & 0 & 0 & \ldots   \\
   0 & a_{2}^{2} & a_{1}^{2} & 0 & 0 & 0 & 0 & \ldots   \\
   0 & a_{3}^{3} & 0 & a_{1}^{3} & 0 & 0 & 0 & \ldots   \\
   0 & a_{4}^{4} & a_{2}^{4} & 0 & a_{1}^{4} & 0 & 0 & \ldots   \\
   0 & a_{5}^{5} & 0 & 0 & 0 & a_{1}^{5} & 0 & \ldots   \\
   0 & a_{6}^{6} & a_{3}^{6} & a_{2}^{6} & 0 & 0 & a_{1}^{6} & \ldots   \\
   \vdots  & \vdots  & \vdots  & \vdots  & \vdots  & \vdots  & \vdots  & \ddots   \\
\end{matrix} \right),  \quad C=\left( \begin{matrix}
   \varphi & 0 & 0 & 0 & 0 & 0 & 0 & \ldots   \\
   0 & c_{1}^{1} & 0 & 0 & 0 & 0 & 0 & \ldots   \\
   0 & c_{2}^{2} & c_{1}^{2} & 0 & 0 & 0 & 0 & \ldots   \\
   0 & c_{3}^{3} & 0 & c_{1}^{3} & 0 & 0 & 0 & \ldots   \\
   0 & c_{4}^{4} & c_{2}^{4} & 0 & c_{1}^{4} & 0 & 0 & \ldots   \\
   0 & c_{5}^{5} & 0 & 0 & 0 & c_{1}^{5} & 0 & \ldots   \\
   0 & c_{6}^{6} & c_{3}^{6} & c_{2}^{6} & 0 & 0 & c_{1}^{6} & \ldots   \\
   \vdots  & \vdots  & \vdots  & \vdots  & \vdots  & \vdots  & \vdots  & \ddots   \\
\end{matrix} \right).$$
It is obvious that ($n>0$)
$$\left[ n,\to  \right]A=\left[ n,\to  \right]\left\langle {{a}^{\left( \ln n \right)}}\left( x \right),x \right\rangle ,$$
$$\left[ n,\to  \right]C=\left[ n,\to  \right]\left\langle x-{{\left( \log \circ a\left( x \right) \right)}^{*}}\circ {{a}^{\left( \ln n \right)}}\left( x \right),x \right\rangle .$$
Since
$$\left( x-{{a}^{*}}\left( x \right)\circ {{a}^{\left( -1 \right)}}\left( x \right) \right)\circ {{a}^{\left( \ln m \right)}}\left( x \right)={{a}^{\left( \ln m \right)}}\left( x \right)-\frac{1}{\ln m}{{\left( {{a}^{\left( \ln m \right)}}\left( x \right) \right)}^{*}},$$
or
$$\left[ {{x}^{n}} \right]\left( x-{{\left( \log \circ a\left( x \right) \right)}^{*}} \right)\circ {{a}^{\left( \ln m \right)}}\left( x \right)=\frac{\ln \left( {m}/{n}\; \right)}{\ln m}\left[ {{x}^{n}} \right]{{a}^{\left( \ln m \right)}}\left( x \right),$$
then
$$\left[ {{x}^{nm}} \right]A{{x}^{m}}\circ \left( x-{{\left( \log \circ a\left( x \right) \right)}^{*}} \right)\circ {{a}^{\left( -\ln m \right)}}\left( x \right)=\left[ {{x}^{nm}} \right]C{{x}^{m}}\circ {{a}^{\left( -\ln m \right)}}\left( x \right)=$$
$$=\left[ {{x}^{n}} \right]\left( x-{{\left( \log \circ a\left( x \right) \right)}^{*}} \right)\circ {{a}^{\left( \ln n \right)}}\left( x \right)=1, \qquad n=1; \qquad=0, \qquad n>1.$$
Thus, up to the element equal $\varphi $, 
$$A=\left\langle x+{{\left( \log \circ b\left( x \right) \right)}^{*}},b\left( x \right) \right\rangle ,  \qquad C=\left\langle x,b\left( x \right) \right\rangle ,$$
$$\left[ {{x}^{n}} \right]{{b}^{\left( \ln m \right)}}\left( x \right)=c_{n}^{mn}=\frac{\ln m}{\ln mn}\left[ {{x}^{n}} \right]{{a}^{\left( \ln mn \right)}}\left( x \right).$$
Denote
$${{\left\langle x,{{a}^{\left( -\beta  \right)}}\left( x \right) \right\rangle }^{-1}}=\left\langle x,{}_{\left( \beta  \right)}{{a}^{\left( \beta  \right)}}\left( x \right) \right\rangle .$$
Then
$$\left[ {{x}^{n}} \right]{}_{\left( \beta  \right)}{{a}^{\left( \beta \ln m \right)}}\left( x \right)=\frac{\beta \ln m}{\beta \ln m+\beta \ln n}\left[ {{x}^{n}} \right]{{a}^{\left( \beta \ln m+\beta \ln n \right)}}\left( x \right).$$
Denote
$$\left[ n,\to  \right]{{\left\langle x|\log \circ {}_{\left( \beta  \right)}a\left( x \right) \right\rangle }_{{{e}^{x}}}}={}_{\left( \beta  \right)}{{\tilde{s}}_{n}}\left( x \right), \qquad_{\left( 0 \right)}{{\tilde{s}}_{n}}\left( x \right)={{\tilde{s}}_{n}}\left( x \right).$$
Then
$$_{\left( \beta  \right)}{{a}^{\left( \varphi  \right)}}\left( x \right)=\sum\limits_{n=1}^{\infty }{\frac{\varphi }{\varphi +\beta \ln n}}\frac{{{{\tilde{s}}}_{n}}\left( \varphi +\beta \ln n \right)}{n!}{{x}^{n}},$$
$$_{\left( \beta  \right)}{{\tilde{s}}_{n}}\left( x \right)=x{{\left( x+\beta \ln n \right)}^{-1}}{{\tilde{s}}_{n}}\left( x+\beta \ln n \right).$$
{\bfseries Example.}
$$\left[ {{x}^{n}} \right]{{\varepsilon }^{\left( \varphi  \right)}}\left( x \right)=\frac{{{\varphi }^{s\left( n \right)}}}{f\left( n \right)},  \qquad\left[ {{x}^{n}} \right]{}_{\left( 1 \right)}{{\varepsilon }^{\left( \varphi  \right)}}\left( x \right)=\frac{\varphi {{\left( \varphi +\ln n \right)}^{s\left( n \right)-1}}}{f\left( n \right)},$$
$$s\left( 1 \right)=0, \qquad f\left( 1 \right)=1, \qquad s\left( n \right)={{m}_{1}}+{{m}_{2}}+...+{{m}_{r}},  \qquad f\left( n \right)={{m}_{1}}!{{m}_{2}}!...{{m}_{r}}!,$$
$n=p_{1}^{{{m}_{1}}}p_{2}^{{{m}_{2}}}...\text{ }p_{r}^{{{m}_{r}}} $ is the canonical decomposition of number $n$. From
$$_{\left( 1 \right)}{{\varepsilon }^{\left( \varphi +\beta  \right)}}\left( x \right)={}_{\left( 1 \right)}{{\varepsilon }^{\left( \varphi  \right)}}\left( x \right)\circ {}_{\left( 1 \right)}{{\varepsilon }^{\left( \beta  \right)}}\left( x \right)$$ 
we obtain analog of the Abel's generalized binomial formula:
$$\left( \varphi +\beta  \right){{\left( \varphi +\beta +\ln n \right)}^{s\left( n \right)-1}}=\sum\limits_{d|n}{{{\left( \begin{matrix}
   n  \\
   d  \\
\end{matrix} \right)}_{f}}}\varphi {{\left( \varphi +\ln d \right)}^{s\left( d \right)-1}}\beta {{\left( \beta +\ln \left( {n}/{d}\; \right) \right)}^{s\left( {n}/{d}\; \right)-1}},$$
where
$${{\left( \begin{matrix}
   n  \\
   d  \\
\end{matrix} \right)}_{f}}=\frac{f\left( n \right)}{f\left( d \right)f\left( {n}/{d}\; \right)};$$ 
or, since 
$$\left[ {{x}^{n}} \right]\left( x+{{\left( \log \circ {}_{\left( 1 \right)}\varepsilon \left( x \right) \right)}^{*}} \right)\circ {}_{\left( 1 \right)}{{\varepsilon }^{\left( \varphi  \right)}}\left( x \right)=\left[ {{x}^{n}} \right]{{\varepsilon }^{\left( \varphi +\ln n \right)}}\left( x \right),$$
then
$${{\left( \varphi +\beta +\ln n \right)}^{s\left( n \right)}}=\sum\limits_{d|n}{{{\left( \begin{matrix}
   n  \\
   d  \\
\end{matrix} \right)}_{f}}}{{\left( \varphi +\ln d \right)}^{s\left( d \right)}}\beta {{\left( \beta +\ln \left( {n}/{d}\; \right) \right)}^{s\left( {n}/{d}\; \right)-1}}.$$
Since
$$_{\left( 1 \right)}{{\varepsilon }^{\left( \varphi  \right)}}\left( x \right)=\left\langle x,{}_{\left( 1 \right)}\varepsilon \left( x \right) \right\rangle {{\varepsilon }^{\left( \varphi  \right)}}\left( x \right),  \qquad{{\varepsilon }^{\left( \varphi  \right)}}\left( x \right)=\left\langle x,{{\varepsilon }^{\left( -1 \right)}}\left( x \right) \right\rangle {}_{\left( 1 \right)}{{\varepsilon }^{\left( \varphi  \right)}}\left( x \right),$$
then
$$\varphi {{\left( \varphi +\ln n \right)}^{s\left( n \right)-1}}=\sum\limits_{d|n}{{{\left( \begin{matrix}
   n  \\
   d  \\
\end{matrix} \right)}_{f}}{{\varphi }^{s\left( d \right)}}\ln d{{\left( \ln n \right)}^{s\left( {n}/{d}\; \right)-1}}},$$
$${{\varphi }^{s\left( n \right)}}=\sum\limits_{d|n}{{{\left( \begin{matrix}
   n  \\
   d  \\
\end{matrix} \right)}_{f}}\varphi {{\left( \varphi +\ln d \right)}^{s\left( d \right)-1}}{{\left( \ln \left( {1}/{d}\; \right) \right)}^{s\left( {n}/{d}\; \right)}}}.$$
When $n={{p}^{m}}$ formulas take the form of Abel's identities [9; p. 92-99]:
$$\left( \varphi +\beta  \right){{\left( \varphi +\beta +ma \right)}^{m-1}}=\sum\limits_{k=0}^{m}{\left( \begin{matrix}
   m  \\
   k  \\
\end{matrix} \right)\varphi {{\left( \varphi +ka \right)}^{k-1}}\beta {{\left( \beta +\left( m-k \right)a \right)}^{m-k-1}}},$$
$${{\left( \varphi +\beta +ma \right)}^{m}}=\sum\limits_{k=0}^{m}{\left( \begin{matrix}
   m  \\
   k  \\
\end{matrix} \right){{\left( \varphi +ka \right)}^{k}}}\beta {{\left( \beta +\left( m-k \right)a \right)}^{m-k-1}},$$
$$\varphi {{\left( \varphi +ma \right)}^{m-1}}=\sum\limits_{k=0}^{m}{\left( \begin{matrix}
   m  \\
   k  \\
\end{matrix} \right)}{{\varphi }^{k}}ka{{\left( ma \right)}^{m-k-1}},$$
$${{\varphi }^{m}}=\sum\limits_{k=0}^{m}{\left( \begin{matrix}
   m  \\
   k  \\
\end{matrix} \right)\varphi {{\left( \varphi +ka \right)}^{k-1}}}{{\left( -ka \right)}^{m-k}},  \qquad a=\ln p.$$

We generalize this example. Let ${{s}_{n}}\left( x \right)$ is the certain binomial sequence. Construct the series $a\left( x \right)$ :
$$\left[ {{x}^{n}} \right]{{a}^{\left( \varphi  \right)}}\left( x \right)=\frac{{{u}_{n}}\left( \varphi  \right)}{f\left( n \right)},
\qquad{{u}_{n}}\left( x \right)={{s}_{{{m}_{1}}}}\left( x \right){{s}_{{{m}_{2}}}}\left( x \right)...{{s}_{{{m}_{r}}}}\left( x \right),  \qquad n=p_{1}^{{{m}_{1}}}p_{2}^{{{m}_{2}}}...\text{ }p_{r}^{{{m}_{r}}}.$$
Then
$$\left[ {{x}^{n}} \right]{}_{\left( \beta  \right)}{{a}^{\left( \varphi  \right)}}\left( x \right)=\frac{\varphi }{\varphi +\beta \ln n}\frac{{{u}_{n}}\left( \varphi +\beta \ln n \right)}{f\left( n \right)},$$
$$_{\left( \beta  \right)}{{a}^{\left( \varphi  \right)}}\left( x \right)=\left\langle x,{}_{\left( \beta  \right)}{{a}^{\left( \beta  \right)}}\left( x \right) \right\rangle {{a}^{\left( \varphi  \right)}}\left( x \right),  \qquad{{a}^{\left( \varphi  \right)}}\left( x \right)=\left\langle x,{{a}^{\left( -\beta  \right)}}\left( x \right) \right\rangle {}_{\left( \beta  \right)}{{a}^{\left( \varphi  \right)}}\left( x \right),$$
$$\frac{\varphi }{\varphi +\beta \ln n}{{u}_{n}}\left( \varphi +\beta \ln n \right)=\sum\limits_{d|n}{{{\left( \begin{matrix}
   n  \\
   d  \\
\end{matrix} \right)}_{f}}{{u}_{d}}\left( \varphi  \right)\frac{\ln d}{\ln n}}{{u}_{{n}/{d}\;}}\left( \beta \ln n \right),$$
$${{u}_{n}}\left( \varphi  \right)=\sum\limits_{d|n}{{{\left( \begin{matrix}
   n  \\
   d  \\
\end{matrix} \right)}_{f}}\frac{\varphi }{\varphi +\beta \ln d}{{u}_{d}}\left( \varphi +\beta \ln d \right)}{{u}_{{n}/{d}\;}}\left( \beta \ln \left( {1}/{d}\; \right) \right).$$
Since ${{u}_{{{p}^{m}}}}\left( x \right)={{s}_{m}}\left( x \right)$, when $n={{p}^{m}}$ formulas take the form of mutually inverse relations for the Lagrange series:
$$\frac{\varphi }{\varphi +ma}{{s}_{m}}\left( \varphi +ma \right)=\sum\limits_{k=0}^{m}{\left( \begin{matrix}
   m  \\
   k  \\
\end{matrix} \right)}{{s}_{k}}\left( \varphi  \right)\frac{k}{m}{{s}_{m-k}}\left( ma \right),$$
$${{s}_{m}}\left( \varphi  \right)=\sum\limits_{k=0}^{m}{\left( \begin{matrix}
   m  \\
   k  \\
\end{matrix} \right)}\frac{\varphi }{\varphi +ka}{{s}_{k}}\left( \varphi +ka \right){{s}_{m-k}}\left( -ka \right),  \qquad a=\beta \ln p.$$

Note the identities for the coefficients ${{n\choose d}_{f}}$, similar to the identities
$$\sum\limits_{k=0}^{n}{\left( \begin{matrix}
   n  \\
   k  \\
\end{matrix} \right)}={{2}^{n}}; \qquad\sum\limits_{k=0}^{n}{\left( \begin{matrix}
   n  \\
   k  \\
\end{matrix} \right)k{{\left( -1 \right)}^{n-k}}}=0, \qquad n\ne 1.$$
Since $\varepsilon \left( x \right)\circ \varepsilon \left( x \right)={{\varepsilon }^{\left( 2 \right)}}\left( x \right)$, ${{\varepsilon }^{*}}\left( x \right)\circ {{\varepsilon }^{\left( -1 \right)}}\left( x \right)={{\left( \log \circ \varepsilon \left( x \right) \right)}^{*}}$, then
$$\sum\limits_{d|n}{{{\left( \begin{matrix}
   n  \\
   d  \\
\end{matrix} \right)}_{f}}}={{2}^{s\left( n \right)}};  \qquad\sum\limits_{d|n}{{{\left( \begin{matrix}
   n  \\
   d  \\
\end{matrix} \right)}_{f}}}\ln d{{\left( -1 \right)}^{s\left( {n}/{d}\; \right)}}=0, \qquad n\ne p.$$

Generalization of the theorem 2 is the formula
$$b\left( x \right)=\left( x-{{\left( \log \circ a\left( x \right) \right)}^{*}} \right)\circ \sum\limits_{n=1}^{\infty }{{{a}^{\left( -\ln n \right)}}\left( {{x}^{n}} \right)\left[ {{x}^{n}} \right]}b\left( x \right)\circ {{a}^{\left( \ln n \right)}}\left( x \right),$$
which follows from
$$\left[ n,\to  \right]{{\left\langle x-{{\left( \log \circ a\left( x \right) \right)}^{*}},{{a}^{\left( -1 \right)}}\left( x \right) \right\rangle }^{-1}}=\left[ n,\to  \right]\left\langle {{a}^{\left( \ln n \right)}}\left( x \right),x \right\rangle .$$

E-mail: {evgeniy\symbol{"5F}burlachenko@list.ru}
\end{document}